\documentclass[11pt]{amsart}
\usepackage{tikz}
\usepackage{amsthm,amssymb}
\usepackage{mathrsfs}
\usepackage{amscd,verbatim}
\usepackage{microtype}
\allowdisplaybreaks[2]
\usepackage{floatflt}
\usepackage{amssymb}
\usepackage{tikz-cd}
\usepackage[all]{xy}
\usepackage[a4paper,margin=27mm,headheight=14pt]{geometry}

\usepackage[colorlinks,linkcolor=blue,citecolor=blue,urlcolor=red]{hyperref}

\renewcommand{\AA}{\mathbb{A}}

\renewcommand{\L}{\mathbb{L}}
\newcommand{\N}{\mathbb{N}}
\newcommand{\PP}{\mathbb{P}}
\newcommand{\Q}{\mathbb{Q}}
\newcommand{\Z}{\mathbb{Z}}

\newcommand{\sI}{\mathcal{I}}

\newcommand{\sO}{\mathcal{O}}

\newcommand*{\xdashrightarrow}[2][]{%
  \mathrel{%
    \mathpalette{\da@xarrow{#1}{#2}{}\da@rightarrow{\,}{}}{}%
  }%
}

\newcommand{\Coh}{\operatorname{Coh}}

\newcommand{\Cone}{\operatorname{Cone}}

\newcommand{\ul}[1]{{\underline{#1}}}

\newcommand{\ulMDM}{\operatorname{\mathbf{\underline{M}DM}}}

\newcommand{\Solid}{\text{Solid}}

\newcommand{\Solidm}{\text{Solid}^m}

\newcommand{\Hom}{\operatorname{Hom}}

\newcommand{\RHom}{\operatorname{R\underline{Hom}}}

\newcommand{\Spec}{\operatorname{Spec}}

\newcommand{\eff}{{\operatorname{eff}}}

\newcommand{\Zar}{{\operatorname{Zar}}}

\newcommand{\id}{{\operatorname{Id}}}

\renewcommand{\lim}{\operatornamewithlimits{\varprojlim}}
\newcommand{\colim}{\operatornamewithlimits{\varinjlim}}

\newcommand{\ol}{\overline}

\renewcommand{\phi}{\varphi}
\renewcommand{\epsilon}{\varepsilon}

\newcommand{\gm}{{\operatorname{gm}}}

\newcommand{\Bl}{{\mathbf{Bl}}}

\newcommand{\bcube}{{\ol{\square}}}

\newcommand{\M}{\mathbf{M}}

\newcommand{\ulM}{\underline{\M}}

\newcommand{\ulPSCH}{\operatorname{\mathbf{\underline{P}Sch}}}

\newcommand{\ulMDMgm}{\ulMDM_{\gm}^{\eff}}

\newcommand{\ulPNor}{\operatorname{\mathbf{\underline{P}Nor}}}
\newcommand{\ulPReg}{\operatorname{\mathbf{\underline{P}Reg}}}
\newcommand{\bbox}{\blacksquare}
\newcommand{\fib}{\operatorname{fib}}
\newcommand{\cofib}{\operatorname{cofib}}

\newcounter{spec}
{\end{list}}%

\newtheorem{lemma}{Lemma}[section]

\newtheorem{Observation}[lemma]{Observation}

\theoremstyle{definition}
\newtheorem{thm}[lemma]{Theorem}

\newtheorem{prop}[lemma]{Proposition}

\theoremstyle{remark}

\newtheorem{remark}[lemma]{Remark}

\numberwithin{equation}{section}

\begin{document}
\title{Solid realization of motives with modulus}
\author{Keiho Matsumoto}

\begin{abstract}
We construct a covariant realization functor $\Solidm(-)$ from the category of motives with modulus to the derived category of solid modules in the sense of Clausen–Scholze. For a smooth modulus pair $(X,D)$, the dual of $\Solidm(X,D)$ recovers the Kelly–Miyazaki Hodge realization of $(X,D)$. Let $X$ be a smooth proper variety over a field, let $D\subset X$ be an effective Cartier divisor whose support has simple normal crossings, and let $U\subset X$ be an open subscheme. We compute $\Solidm\!\bigl(\Cone(\ulM(U,D|_{U})\to \ulM(X,D))\bigr)$ explicitly in terms of the formal completion along $X\setminus U$ of a model of $X$ of finite type over $\Z$.

\end{abstract}

\maketitle

\tableofcontents

\section{Introduction}

For a discrete Noetherian commutative ring $A$, let $\ulPNor(A)$ be the full subcategory of $\ulPSCH(A)$ consisting of modulus pairs with a normal total space. Consider the pre-analytic ring naturally associated with $A$, $A_{\bbox}:=(A,A)_{\bbox}$, where the underlying condensed ring is $A$. For an extremally disconnected set $S=\lim S_i$, the corresponding condensed module is $\displaystyle\colim_{A'\to A|A':\text{f.g.alg}/\mathbb{Z}} \lim A'_{\bbox}[S_i]$. We denote the derived category of solid $A_{\bbox}$-modules by $D(A_{\bbox})$. In this case, we introduce a {covariant} functor:
\[
\text{Solid}:\ulPNor(A) \to D(A_{\bbox})
\]
which maps a modulus pair $(X, D) \in \ulPNor(A)$ to a solid $A_{\bbox}$-complex:
\[
f_!(f^!A \otimes^{\mathbb{L}}_{\sO_{X,\bbox}}\sO_{X}(-D + |D|)),
\]
where $f$ denotes the structure map $X \to \Spec A$. Specifically, we will show that this functor has good properties, such as $\bcube$-invariance, smooth blow-up invariance and symmetric monoidality. We note that $\text{Solid}(X,D)$ can be viewed as the $R_{{\sO}}(``M(X,D)^{pre-dual}")$, where $R_\sO$ is the Kelly-Miyazaki's realization functor $\ulMDM^\eff_\gm(A) \to D(A)^{\text{op}}~|~\ulM(X,D) \mapsto R\Gamma_{\Zar}(X,\sO_X(D-|D|))$ constructed in \cite{kelly2023hodge}, and $``\ulM(X,D)^{pre-dual}"$ is the ``motive whose dual is $\ulM(X,D)$" in the suitable category.

\begin{thm}\label{main}
     Let $A$ be a discrete Noetherian commutative ring. Let $A_{\bbox}$ be the analytic ring $(A, A)_{\bbox}$, and let $D(A_{\bbox})$ be the derived category of solid $A_{\bbox}$-modules. 
    \begin{itemize}
        \item[(1)] Let $h: (Y, C) \to (X, D)$ be an ambient morphism over $(\Spec A, \emptyset)$. Then, there exists a morphism in $D(A_{\bbox})$:
        \[
        h_+: \text{Solid}(Y, C) \to \text{Solid}(X, D).
        \]
        Furthermore, for a sequence of ambient morphisms $(Z, B) \overset{k}{\to} (Y, C) \overset{h}{\to} (X, D)$, the functoriality $(h \circ k)_+ = h_+ \circ k_+$ holds.
        \item[(2)] $\bcube$-invariance holds. That is, the natural morphism:
        \[
        \Solid((X, D) \otimes \bcube) \to \Solid(X, D)
        \]
        becomes an isomorphism in $D(A_{\bbox})$.
        \item[(3)] Assume $X$ is regular, $D$ is strictly normal crossing and $X \to \Spec A$ is locally complete intersection. Given an integral regular scheme $Z$ and a regular closed embedding $Z \hookrightarrow X$ with $|Z| \subset |D|$, the invariance under blow-ups holds. That is, the natural morphism:
        \[
        \Solid(\text{Bl}_Z X, p^* D) \to \Solid(X, D)
        \]
        becomes an isomorphism. Here, $p: \text{Bl}_Z X \to X$ denotes the blow-up of $X$ along $Z$.
        \item[(4)] Assume that $X$ is proper smooth over $A$. Then, $\Solid(X, D)$ is a discrete solid $A_{\bbox}$-module complex (that is a usual $A$- module complex), and its $i$-th cohomology group is:
        \[
        H^{i + \dim X/A}(X, \Omega_{X/A}^{\dim X/A}(-D + |D|)).
        \]
        When $A$ is a field, by Serre duality, this coincides with:
        \[
        H^{-i}(X, \sO(D - |D|))^{\vee}.
        \]
        \item[(5)] Assume $X$ is complete intersection over $A$. The dual of $\Solid(X, D)$ over $A$:
        \[
        R\ul{\Hom}_{A_\bbox}(\Solid(X, D), A)
        \]
        is a discrete solid $A_{\bbox}$-module, and it is quasi-isomorphic to 
        \[
        R\Gamma_{\Zar}(X,\sO_X(D-|D|)).
        \]
        \item[(6)] Let $(X,D)$ and $(Y,C)$ be modulus pairs. We assume $X$ and $Y$ are smooth over $A$, then there is an isomorphism
        \[
    \Solid\bigl((X,D)\otimes(Y,C)\bigr)  \simeq\Solid(X,D) \otimes_{A_{\bbox}}^\L \Solid(Y,C) 
        \]
        in $D(A_{\bbox})$.
        \item[(7)] For a $\underline{MV}$-square, which is “elementary Nisnevich square of modulus pairs” (see \cite[Definition 4.1.1]{miyazaki2020nisnevich}),
        \[
        \xymatrix{
        (V,j^*C) \ar[r]^-{j} \ar[d]_{g} & (Y,C) \ar[d]^{f}\\
        (U,i^*D) \ar[r]^-{i} & (X,D)
        }
        \]
        there exists a distinguished triangle 
        \[
        \Solid(V,j^*C) \overset{(j_+,g_+)}{\to} \Solid(Y,C) \oplus \Solid(U,i^*D) \overset{i_+-f_+}{\to} \Solid(X,D) \to 
        \]
        \item[(8)] Assume $A$ is a field of characteristic $0$. There is a covariant functor 
        \[
        \ulMDM_\gm^\eff(A,\Q) \to D(A_{\bbox}) 
        \]
        which sends $\ulM(X,D)_\Q$ to $\Solid(X,D)$.
    \end{itemize}
\end{thm}

\begin{remark}
    We provide some remarks on Theorem (5) above. In the case where $X$ is not proper, a $A$-module complex $P$ satisfying $\RHom_A(P,A)\simeq R\Gamma_{\Zar}(X,\sO_X(D-|D|))$ may not exist. For example, when $X = \AA_A^1$, it is easy to see that there does not exist a $A$-module complex $P$ such that $\RHom_A(P,A)\simeq R\Gamma_\Zar(\AA^1,\sO_{\AA^1})$ holds. Therefore, Theorem (5) above is a feature unique to the theory of solid modules, which cannot be considered within the framework of discrete theories.
\end{remark}

In Section~3, for a regular scheme $X$ over $A$, a simple normal crossings Cartier divisor $D$ on $X$, and each $m \in \mathbb{Z}_{\geq 0}$, 
we construct a solid module 
\[
\Solidm(X,D) \in D(A_\bbox).
\] 
For $m=0$, this recovers the previously defined object $\Solid(X,D)$. 
The construction is characterized by the property that its dual $ R\ul{\Hom}_{A_\bbox}(\Solidm(X, D), A)$ is discrete, and that the $i$-th cohomology group of this dual coincides with Kelly–Miyazaki’s Hodge cohomology with modulus:
\[
H^i\bigl((\Solidm(X,D))^\vee\bigr)\;\simeq\; H_\Zar^i(X,\Omega^m(\log |D|)(D-|D|)).
\]
We also prove that if the base ring $A$ is a field of characteristic $0$, then $\Solidm$ can be extended to the functor
\[
\Solidm:\ulMDMgm(A,\Q) \to D(A_{\bbox}).
\]

In Section~4, we give an explicit description of $\Solidm(X,D)$ for a smooth variety $X$ over a field $A$. The main input is Proposition~\ref{renprosolid}. It describes $j_!j^*F$ for an open immersion $j$ of schemes of finite type over $\Z$ and a perfect complex $F$ in terms of the formal completion along the closed complement. For a principal open subset defined by a non-zero-divisor, this is Ren's pro-solid comparison theorem \cite{prosolid}. As an application, for a proper smooth variety $X$ over a field $A$, an effective Cartier divisor $D\subset X$ such that $|D|$ is a simple normal crossings divisor, and an open immersion $j:U\hookrightarrow X$ (with closed complement $i:Z:=X\setminus U\hookrightarrow X$), we compute
\[
\Solidm\!\Bigl(\Cone\bigl(\ulM(U,D|_{U})\longrightarrow \ulM(X,D)\bigr)\Bigr)
\]
(Theorem~\ref{solidboundarycone}). In particular, we show that the resulting object admits a geometric description in terms of the formal completion along $Z$ of a model of $X$ of finite type over $\Z$. This is a feature not shared by the Kelly–Miyazaki Hodge realization. We note that the completion has to be taken on such a model, and not on $X$ itself (Remark~\ref{rem:naive}).

\section*{Acknowledgements}
The author would like to thank Hiroyasu Miyazaki for helpful discussions. The author is grateful to Shane Kelly for helpful comments, in particular for his comments on the revision of Section~4. The author also thanks Fei Ren for pointing out an error in Section~4 of the first version of this paper.

\section*{Use of AI}
Based on comments from Shane Kelly, the author revised Section~4 with the assistance of the generative AI tools ChatGPT (GPT-6.1 Astra, OpenAI) and Claude (Claude Opus 5.5, Anthropic). These tools were also used for improving the English of Section~4. The author has checked all the proofs and takes full responsibility for the contents of this paper.


\section{Proof of Theorem~\ref{main}}

\subsection{Review of Exceptional Functors $f^!$, $f_!$} Let $A$ be a discrete commutative ring. Consider the pre-analytic ring naturally associated with $A$, $A_{\bbox}:=(A,A)_{\bbox}$, where the underlying condensed ring is $A$. For an extremally disconnected set $S=\lim S_i$, the corresponding condensed module is $\displaystyle\colim_{A'\to A|A':\text{f.g.alg}/\mathbb{Z}} \lim A'_{\bbox}[S_i]$. We denote the derived category of solid $A_{\bbox}$-modules by $D(A_{\bbox})$. More generally, for a separated scheme $X$ over a ring $A$, we can construct the derived category $D(\sO_{X,\bbox})$ of solid $\sO_{X,\bbox}$-modules (\cite[section 10]{cond}). The most basic example of a solid $\sO_{X,\bbox}$-module is a coherent $\sO_X$-module $\mathcal{F}$. We refer to a solid $\sO_{X,\bbox}$-module obtained from such a coherent $\sO_X$-module as a discrete solid $\sO_{X,\bbox}$-module. Furthermore, there exists a fully faithful functor
\[
D^b(\text{Coh}(X)) \to D(\sO_{X,\bbox})
\]
Clausen and Scholze proved the following theorem:

\begin{thm}[{\cite[affine case in section 8, general case in section 11]{cond}}]\label{CSex}
    Let $R$ be a commutative Noetherian ring, and let $f:Y\to X$ be a morphism between separated, finite type schemes over $R$. Then, the functors $f_!:D(\sO_{Y,\bbox}) \to D(\sO_{X,\bbox})$ and $f^!:D(\sO_{X,\bbox}) \to D(\sO_{Y,\bbox})$ exist, satisfying the following properties:
    \begin{itemize}
        \item[(1)] (Verdier duality) $f_!$ is left adjoint to $f^!$. In particular, for any $M\in D(\sO_{X,\bbox})$ and $N\in D(\sO_{Y,\bbox})$, there is an isomorphism in $D(\sO_{X,\bbox})$:
        \[
        f_*R\ul{\Hom}_{\sO_Y}(N,f^!M) \simeq R\ul{\Hom}_{\sO_X}(f_!N,M)
        \]
        \item[(2)] (Poincaré duality) $f^!\sO_X$ is a discrete solid module. Moreover, if $f$ has finite Tor dimension (i.e., $f_*\sO_Y[0]$ has finite Tor dimension in $D(\text{Coh}(X))$), then for any discrete complex $M \in D(\sO_{X,\bbox})$,
        \[
        f^!M \simeq f^!\sO_X\otimes^{\mathbb{L}}_{\sO_{Y}}  f^*M 
        \]
        holds. Additionally, by the symmetry of derived tensor product and its commutativity with derived pullback, for any discrete complexes $M,N \in D(\sO_{X,\bbox})$,
        \[
        f^!(N\otimes M) \simeq f^! \sO_{X}\otimes^{\mathbb{L}}_{\sO_{Y}} f^*N \otimes^{\mathbb{L}}_{\sO_{Y}}  f^*M \simeq f^!N   \otimes^{\mathbb{L}}_{\sO_{Y}}  f^*M
        \]
        holds.
        \item[(3)] (Projection formula) For any $M\in D(\sO_{X,\bbox})$ and $N\in D(\sO_{Y,\bbox})$, there is an isomorphism in $D(\sO_{X,\bbox})$:
        \[
        f_!(N\otimes^{\mathbb{L}}_{\sO_{Y,\bbox}}f^*M) \simeq f_!N\otimes^{\mathbb{L}}_{\sO_{X,\bbox}}M
        \]
        \item[(4)] Suppose $f$ is smooth. If we denote the relative dimension of $f:Y\to X$ by $n$, then there is an isomorphism:
        \[
        f^!\sO_X \simeq \Omega^{n}_{Y/X}[n].
        \]
        \item[(5)] If $f$ is proper, $f_!$ coincides with $f_*$. In particular, for discrete bounded complexes, $f_!$ coincides with $Rf_*$. That is, the following diagram commutes:
        \[
        \xymatrix{
D^b(\text{Coh}(B)) \ar[d]_{Rf_*} \ar[r]& D(B_{\bbox}) \ar[d]^{f_!=f_*} \\
D^b(\text{Coh}(A))  \ar[r]& D(A_{\bbox})
}
        \]
        \item[(6)] If $f$ is a complete intersection, $f^!\sO_X$ is locally a line bundle concentrated in some degree.
        \item[(7)] Given a morphism $g:S\to X$ between finite type, separated schemes, consider the following Cartesian diagram:
        \[
        \xymatrix{
        S\times_X Y \ar[r]^-{g'} \ar[d]_{f'} & Y \ar[d]^-{f}\\
        S \ar[r]_{g} & X
        }
        \]
        If $g$ is flat, then there is a natural equivalence:
        \[
         g^* \circ f_! \simeq f'_!\circ {g'}^*:D(\sO_{Y,\bbox}) \to D(\sO_{S,\bbox})
        \]
        \item[(8)] Let $X$ be an integral regular scheme and $Y$ a Cartier divisor on $X$. Then the following isomorphism holds:
        \[
        f^! \sO_X = \sO_Y(Y)[-1].
        \]
        Here, $\sO_Y(Y)$ represents a locally free sheaf $f^*\sO_X(Y)$ on $Y$.
        \end{itemize}  
    \end{thm}

   \begin{remark}\label{remark}
        \begin{itemize}
            \item[(1)] Except for (4), (6), and (8), the above theorem is extracted from Clausen-Scholze's theorem (\cite[section 11]{cond}) stating that the correspondence $X \mapsto D(\sO_{X,\bbox})$ forms a 6-functor formalism, which is used in this note. (4) is proven in \cite[Thm 11.6]{cond}. (6) is proven in \cite[Thm 8.13 (iv)]{cond}. (8) is a special case of \cite[Prop 11.5]{cond}.


            \item[(2)] Although \cite{cond} assumes that $R$ is of finite type over $\Z$, this restriction can be removed by a standard Noetherian-approximation argument, so the theorem remains valid for arbitrary commutative Noetherian base rings; see \cite{analyticstack} for details.

            \item[(3)] In this section we will prove theorem~\ref{main} by explicit computations. In Section~3 we construct a more general solid functor $\Solidm(-)$ and establish an analogous theorem; the argument there is more conceptual and relies on the Kelly–Miyazaki Hodge realization \cite{kelly2023hodgecohomologyramificationfiltration}, yielding a streamlined proof.
        \end{itemize}
    \end{remark}

    \subsection{Proof of the Theorem~\ref{main}} 
    
    \begin{Observation}[Proof of Theorem \ref{main} (1)]
     Let $h:(Y,C)\to (X,D)$ be an ambient morphism between modulus pairs where the total space is normal over $(\Spec A,\emptyset)$. Let $g$ be the structure morphism of $Y$ over $A$, and $f$ be the structure morphism of $X$ over $A$. Then, there exists the following natural inclusion:
     \begin{equation}\label{sub}
     \sO_Y(-C+|C|) \hookrightarrow h^*\sO_X(-D+|D|).
     \end{equation}
     Now, by Theorem~\ref{CSex} (3), we obtain the following isomorphism:
     \begin{equation}\label{counit}
     h_!\Bigl(h^!f^!A\otimes_{\sO_{Y,\bbox}}^{\L} \sO_Y(-h^*D+h^*|D|)\Bigr) \simeq h_!h^!f^!A\otimes ^{\mathbb{L}}_{\sO_{X,\bbox}} \sO_X(-D+|D|).
     \end{equation}
Here, from the existence of the counit $h_!h^! \Longrightarrow \text{id}$, there exists a morphism from the right side of \eqref{counit}:
\begin{equation}\label{counit2}
    h_!h^!f^!A\otimes ^{\mathbb{L}}_{\sO_{X},\bbox} \sO_X(-D+|D|) \to f^!A\otimes_{\sO_{X,\bbox}}^{\L} \sO_X(-D+|D|)
\end{equation}
Combining the morphisms \eqref{sub}, \eqref{counit}, and \eqref{counit2}, we obtain the following morphism:
\begin{eqnarray}
    h_!\Bigl(g^!A\otimes_{\sO_{Y,\bbox}}^{\L} \sO_Y(-C+|C|)\Bigr) & \overset{\eqref{sub}}{\to} &  h_!\Bigl(g^!A\otimes_{\sO_{Y,\bbox}}^{\L} \sO_Y(-h^*D+h^*|D|)\Bigr) \nonumber\\
     & \overset{\eqref{counit}}{\simeq} & h_!h^!f^!A\otimes ^{\mathbb{L}}_{\sO_{X,\bbox}} \sO_X(-D+|D|) \nonumber\\
      &\overset{\eqref{counit2}}{\to}&f^!A\otimes_{\sO_{X,\bbox}}^{\L} \sO_X(-D+|D|) \label{hplus}
\end{eqnarray}
Note that $h^!f^! A = g^!A$ (\cite[Thm 8.13 (iv)]{cond}). By applying $f_!$, we obtain the desired morphism:
\[
h_+:g_!\Bigl(g^!A\otimes_{\sO_{Y,\bbox}}^{\L} \sO_Y(-C+|C|)\Bigr) \to f_!\Bigl( f^!A\otimes_{\sO_{X,\bbox}}^{\L} \sO_X(-D+|D|) \Bigr)
\]
Next, we show the functoriality of $\Solid(-)$.

Consider the following commutative diagram of normal modulus pairs:
\[
\xymatrix{
(Z,B) \ar@/^18pt/[rr]^-{r} \ar[r]^{k}\ar[dr]_{j} & (Y,C) \ar[d]^{g} \ar[r]^{h} & (X,D) \ar[ld]^{f} \\
&\Spec A &}
\]
In this case, there exists a commutative diagram formed by natural inclusions:
\[
\xymatrix{
r_!(j^!A \otimes_{\sO_{Z,\bbox}}^{\L} \sO_Z(-B+|B|)) \ar[rd]_-{(3)}\ar[r]^-{(1)} & r_!(j^!A \otimes_{\sO_{Z,\bbox}}^{\L} k^*\sO_Y(-C+|C|)) \ar[d]^-{(2)}\\
& r_!(j^!A \otimes_{\sO_{Z,\bbox}}^{\L} r^*\sO_X(-D+|D|))
}
\]
Here, using the projection formula and counit, we obtain the following diagram:\footnotesize
\[
\resizebox{\textwidth}{!}{$\displaystyle
\xymatrix{
 r_!(j^!A \otimes_{\sO_{Z,\bbox}}^{\L} k^*\sO_Y(-C+|C|)) \ar[d]^-{(2)} \ar[r]^-{(4)}_-{\simeq} & h_!\bigl(k_!j^!A  \otimes_{\sO_{Y,\bbox}}^{\L} \sO_Y(-C+|C|)  \bigr) \ar[d]^-{} \ar[r]^-{(5)} & h_!\bigl(g^!A \otimes_{\sO_{Y,\bbox}}^{\L} \sO_Y(-C+|C|)\bigr) \ar[d]^-{(6)}\\
  r_!(j^!A \otimes_{\sO_{Z,\bbox}}^{\L} r^*\sO_X(-D+|D|)) \ar[d]_-{(9)}^-{\simeq}\ar[r]_-{(7)}^-{\simeq}& h_!\bigl(k_!j^!A  \otimes_{\sO_{Y,\bbox}}^{\L} h^*\sO_X(-D+|D|)  \bigr) \ar[r]^-{(8)}
& h_!\bigl(g^!A \otimes_{\sO_{Y,\bbox}}^{\L} h^*\sO_X(-D+|D|)\bigr) \ar[d]^-{(10)}_-{\simeq}\\
r_!r^!f^!A \otimes_{\sO_{X,\bbox}}^{\L} \sO_X(-D+|D|) \ar[r]_-{(11)} &  f^!A \otimes_{\sO_{X,\bbox}}^{\L} \sO_X(-D+|D|)& h_!g^!A \otimes_{\sO_{X,\bbox}}^{\L} \sO_X(-D+|D|) \ar[l]^-{(12)}
}
$}
\]\normalsize
Here, the vertical arrows on the upper side of the diagram come from inclusions, while (5), (8), (11), and (12) come from the counit. Moreover, (4), (7), (9), and (10) are obtained from the projection formula. Hence the top rectangle commutes.

It remains to prove the commutativity of the bottom rectangle. Since the assertion is Zariski–local on $X$, we may replace $\sO(-D+|D|)$ by $\sO_X$. In this case, the commutativity follows from the natural transformation.
$$
r_!r^!\Longrightarrow \mathrm{id}
$$
given by the counit of $r_!\dashv r^!$, together with the composite of counits for $k_!\dashv k^!$ and $g_!\dashv g^!$:
$$
h_!k_!k^!h^!\Longrightarrow h_!h^!\Longrightarrow \mathrm{id}.
$$
This establishes the commutativity of the entire diagram.

By definition, we have
$$
r_+ \;=\; f_!\Bigl((11)\circ(9)\circ(3)\Bigr),\qquad
k_+ \;=\; f_!\Bigl((5)\circ(4)\circ(1)\Bigr),\qquad
h_+ \;=\; f_!\Bigl((12)\circ(10)\circ(6)\Bigr).
$$
The commutativity of the above diagram then proves the claim.
\end{Observation}

\begin{Observation}[Proof of Theorem~\ref{main} (2)]
For a modulus pair $(X,D) \in \ulPReg(A)$, we show that the following natural morphism:
\[
p_+:\text{Solid}((X,D)\otimes \bcube) \to \text{Solid}(X,D)
\]
is an isomorphism. Here, $p:X\times\PP^1 \to X$. This calculation is interesting and essentially reduces to the projection formula that for any $M\in D(\text{Coh}(X))$, there exists the following isomorphism:
\begin{equation}\label{PP1}
Rp_*(\Omega^1_{X\times \PP^1/X}[1]\otimes_{\sO_{X\times\PP^1}}p^*M ) \simeq M.
\end{equation}
First, recalling the construction of $h_+$, the claim follows by proving that each step \eqref{sub}, \eqref{counit}, and \eqref{counit2} of \eqref{hplus} is an isomorphism. The isomorphism of \eqref{sub} follows from the coincidence of the Weil divisors:
\begin{eqnarray*}
- (X\times\{\infty\} + D\times \PP^1) + |X\times\{\infty\} + D\times \PP^1| &= & -D\times \PP^1 +|D|\times \PP^1 \\
&=& p^*(-D +|D|)
\end{eqnarray*}
\eqref{counit} is always an isomorphism. We show the isomorphism of \eqref{counit2}. By Theorem~\ref{CSex} (2), (4), (5), we obtain the following sequence of isomorphisms:
\begin{align*}
   &p_! p^!(f^!A \otimes_{\sO_{X,\bbox}}^{\L} \sO_X(-D+|D|))\\
   &\quad\overset{(2)}{=} p_!\Bigl( p^*(f^!A \otimes_{\sO_{X,\bbox}}^{\L} \sO_X(-D+|D|)) \otimes^{\mathbb{L}}_{\sO_Y} p^!\sO_X \Bigr)\\
    &\quad\overset{(4),(5)}{\simeq} Rp_*\Bigl( p^*(f^!A \otimes_{\sO_{X,\bbox}}^{\L} \sO_X(-D+|D|)) \otimes_{\sO_{Y,\bbox}}^{\L} \Omega^1_{X\times\PP^1/X}[1] \Bigr) \\
    &\quad\overset{\eqref{PP1}}{\simeq} f^!A \otimes_{\sO_{X,\bbox}}^{\L} \sO_X(-D+|D|).
\end{align*}
This indeed shows the isomorphism of \eqref{counit2}.
\end{Observation}

\begin{Observation}[Proof of Theorem~\ref{main} (3)] 
Assume $X$ is regular, $D\hookrightarrow X$ is a regular embedding, and $X \to \Spec A$ is a locally complete intersection. Given an integral regular scheme $Z$ and a regular closed embedding $Z \hookrightarrow X$ with $|Z| \subset |D|$, then $\Bl_Z X$ is also regular. Denote by $f$ the morphism $X \to \Spec A$, by $g$ the morphism $\Bl_ZX \to \Spec A$, and by $p$ the morphism $\Bl_Z X \to X$. Then, we show the following isomorphism in $D(\sO_{X,\bbox})$:
\begin{equation}\label{prebl}
p_!(g^!A \otimes_{\sO_{\Bl_ZX,\bbox}}^{\L} \sO_{\Bl_ZX}(-p^*D+|p^*D|)) \simeq f^!A\otimes_{\sO_{X,\bbox}}^{\L} \sO_X(-D+|D|)
\end{equation}
The morphism from the left-hand side to the right-hand side is given by \eqref{hplus}. By applying $f_!$, we derive the claim of Theorem~\ref{main} (3). Let $E$ denote the exceptional divisor of $p:\Bl_Z X \to X$, and let $c$ denote the codimension of the regular embedding $Z\hookrightarrow X$. We have the following commutative diagram:
\[
\xymatrix{
E\ar[r]^-{i'}\ar[d]_{p'}& \Bl_Z X \ar[d]^{p} \ar[dr]^{g}&\\
Z \ar[r]_{i}& X \ar[r]_{f}& \Spec A
}
\]
First, we prove the following Proposition.
\begin{prop}\label{discrepancy}
    $p^!\sO_X=\sO_{\Bl_ZX}((c-1)E)$.
\end{prop}
\begin{proof}
    The claim is local on $X$, so there exists a commutative ring $R$ and a regular sequence $(e_1,e_2,...,e_c)$ such that $X=\Spec R$ and $Z= \Spec R/(e_1,e_2,...,e_c)$. In this case, $i^! \sO_X=\sO_Z[-c]$ (see \cite[Prop 11.5]{cond}). Now, since $p$ is an isomorphism over $X\setminus Z$, we have $p^!\sO_{X}|_{\Bl_ZX\setminus E}\simeq \sO_{\Bl_ZX\setminus E}$, and since $Z$ is integral, $E$ is also integral, hence there exists an integer $n$ such that $p^!\sO_{X}\simeq\sO_{\Bl_Z X}(nE)$. Using Poincaré duality:
    \[
    {i'}^!p^!\sO_X={i'}^!\sO_{\Bl_ZX} \otimes_{\sO_{E,\bbox}}^\L {i'}^*p^!\sO_X
    \]
    where \cite[Prop 11.5]{cond} gives ${i'}^!\sO_{\Bl_ZX} \simeq {i'}^*\sO_{\Bl_ZX}(E)[-1]$, and since $p^!\sO_{X}\simeq\sO_{\Bl_Z X}(nE)$, we have:
    \[
     {i'}^!p^!\sO_X= {i'}^*\Bigl(\sO_{\Bl_Z X}((n+1)E)[-1] \Bigr)=\sO_E(-(n+1))[-1].
    \]
    On the other hand, since $ {i'}^!p^!\sO_X= {p'}^!i^!\sO_X$, using Poincaré duality:
    \[
     {i'}^!p^!\sO_X= {p'}^!i^!\sO_X = {p'}^!\sO_Z \otimes_{\sO_{E,\bbox}}^\L {p'}^*i^!\sO_X.
    \]
    Since $p'$ is a smooth morphism of relative dimension $c-1$, ${p'}^!\sO_Z = \Omega^{c-1}_{E/Z}[c-1]$, and since $i^! \sO_X=\sO_Z[-c]$, we have:
    \[
       {i'}^!p^!\sO_X= \Omega^{c-1}_{E/Z}[-1] = \sO_E(-c)[-1].
    \]
    Therefore, $-c=-n-1$, hence $n=c-1$.
\end{proof}
    Now, let $\{D_1,D_2,...,D_j\}$ be the set of irreducible components of $|D|$ containing $|Z|$, and let $\{C_1,C_2,...,C_l\}$ be the set of irreducible components of $|D|$ not containing $|Z|$. Then, for $1\leq a \leq j$, $p^*D_a=p^{-1}D_a+E$, and for $1\leq b \leq l$, $p^*C_b=p^{-1}C_b$, where $p^{-1}(-)$ denotes the strict transform of the divisor along $p$. Since $|D|$ is a strictly normal crossing divisor, $j\leq c$. Thus, we have the following calculation:
    \[
    p^*(|D|)-|p^*D|= (j-1)E.
    \]
Now, we perform the following calculation:
\begin{eqnarray*}
&&p_!(g^!A \otimes_{\sO_{\Bl_ZX,\bbox}}^{\L} \sO_{\Bl_ZX}(-p^*D+|p^*D|))\\
&\overset{\text{Poincaré}}{=}&p_!(p^!\sO_X  \otimes_{\sO_{\Bl_ZX,\bbox}}^{\L} p^* f^!A \otimes_{\sO_{\Bl_ZX,\bbox}}^{\L} \sO_{\Bl_ZX}(-p^*D+p^*|D|-(j-1)E))\\
&\overset{\text{proj.form}}{=}& p_!(p^!\sO_X\otimes_{\sO_{\Bl_ZX,\bbox}}^{\L} \sO_{\Bl_ZX}(-(j-1)E)) \otimes_{\sO_{X,\bbox}}^\L f^!A  \otimes_{\sO_{X,\bbox}}^\L \sO_X(-D+|D|)\\
&\overset{\text{Prop}~\ref{discrepancy}}{=}& p_!\sO_{\Bl_ZX}((c-j)E)  \otimes_{\sO_{X,\bbox}}^\L f^!A  \otimes_{\sO_{X,\bbox}}^\L \sO_X(-D+|D|)
\end{eqnarray*}
Thus, we show that $p_!\sO_{\Bl_ZX}((c-j)E) = \sO_X$. Since $\sO_{\Bl_ZX}((c-j)E)$ is a discrete bounded complex and $p$ is projective, we have $p_!\sO_{\Bl_ZX}((c-j)E)=Rp_*\sO_{\Bl_ZX}((c-j)E)$. Since $X$ is normal, $R^0p_*\sO_{\Bl_ZX}((c-j)E)=\sO_X$. By using Grothendieck formal function theorem and the vanishing of coherent cohomology groups for projective spaces, it is easy to see that $R^np_*\sO_{\Bl_ZX}((c-j)E)=0$ for $n>0$.
\end{Observation}

Theorem~\ref{main} (4) can be directly proven using Theorem~\ref{CSex} (4), (5).

\begin{Observation}[Proof of Theorem~\ref{main} (5)]
Assume $f:X\to \Spec A $ is a complete intersection. By Theorem~\ref{CSex} (1),
\[
R\ul{\Hom}_{A_\bbox}(\Solid(X, D), A) = f_*\ul{\Hom}_{\sO_{X,\bbox}}(f^!A\otimes_{\sO_{X,\bbox}}\sO_X(-D+|D|),f^!A).
\]
By Theorem~\ref{CSex} (6), $f^!A$ is invertible, we have:
\[
R\ul{\Hom}_{A_\bbox}(\Solid(X, D), A) = f_*\ul{\Hom}_{\sO_{X,\bbox}}(\sO_X(-D+|D|),\sO_X)=f_*\sO_{X}(D-|D|).
\]
Since $\sO_{X}(D-|D|)$ is discrete, $f_*\sO_{X}(D-|D|)=Rf_*\sO_{X}(D-|D|)$.
\end{Observation}

\begin{Observation}[Proof of Theorem~\ref{main} (6)]\label{observation6} Consider the following Cartesian diagram:
\[
\xymatrix{
Z:=X\times_{\Spec A} Y \ar[d]_-{f'} \ar[r]^-{g'} & X\ar[d]^{f} \\
Y \ar[r]_{g} & \Spec A
}
\]
We perform the following calculation, using $-\otimes -$ for all tensor products:
\begin{align*}
   &\Solid(X,D) \otimes \Solid(Y,C)\\
   &\quad= f_!(f^!A\otimes \sO_X(-D+|D|)) \otimes \Solid(Y,C) \\
    &\quad\overset{\text{proj.form}}{=} f_!\Bigl(f^!A\otimes \sO_X(-D+|D|) \otimes f^*\Solid(Y,C) \Bigr)\\
    &\quad= f_!\Bigl(f^!A\otimes \sO_X(-D+|D|) \otimes f^*g_!(g^!A\otimes \sO_Y(-C+|C|)) \Bigr)\\
    &\quad\overset{\text{base.change}}{=} f_!\Bigl(f^!A\otimes \sO_X(-D+|D|) \otimes g'_!{f'}^*(g^!A\otimes \sO_Y(-C+|C|)) \Bigr) \\
    &\quad\overset{\text{proj.form}}{=} f_!g'_!\Bigl({g'}^*(f^!A\otimes \sO_X(-D+|D|)) \otimes {f'}^*(g^!A\otimes \sO_Y(-C+|C|)) \Bigr)\\
    &\quad= f_!g'_!\Bigl({g'}^*\Omega^n_{X/A} \otimes {f'}^*\Omega^l_{Y/A} \otimes \sO_{X\times Y}(-B+|B|) \Bigr)[n+l]
\end{align*}
where $B$ is the Cartier divisor $X\times C + D\times Y$ on $X\times_{\Spec A} Y$, $n=\dim X/A$, and $l=\dim Y/A$. Now, since ${g'}^*\Omega^n_{X/A} \otimes {f'}^*\Omega^l_{Y/A} = \Omega^{n+l}_{X\times_{\Spec A} Y/A}$, the claim is proven.
\end{Observation}

\begin{Observation}[Proof of Theorem~\ref{main} (7)] 
Theorem~\ref{main} (7) immediately follows from the following proposition. 

\begin{prop}\label{solidNisdescent}
    Given an elementary Nisnevich square
    \[
    \xymatrix{
    V \ar[r]^{j} \ar[d]_{g} & Y \ar[d]^{f} \\
    U \ar[r]_{i} & X
    }
    \]
    and $M \in D(\sO_{X,\bbox})$. Then $M$ is isomorphic to the cone of
    \[
    f_!j_!j^!f^! M \overset{(c_1,c_2)}{\to} f_!f^!M \oplus i_!i^!M,
    \]
    where 
\[
c_1:f_!j_!j^!f^! M  \to f_!f^!M 
\]
is induced by counit of adjunction $j_! \dashv j^!$, and 
\[
c_2:f_!j_!j^!f^! M = i_!g_!g^!i^!M \to  i_!i^!M
\]
is also induced by counit of adjunction $g_! \dashv g^!$.
\end{prop}
\begin{proof}[Proof of Proposition \ref{solidNisdescent}]
Let us denote $C=\Cone\bigl(f_!j_!j^!f^! M \overset{(c_1,c_2)}{\to} f_!f^!M \oplus i_!i^!M\bigr)$. It is suffice to show that for any $P \in D(\sO_{X,\bbox})$ the natural map
    \begin{equation}\label{descentNismap}
 \ul{\Hom}_{\sO_{X,\bbox}}(M,P) \to  \ul{\Hom}_{\sO_{X,\bbox}}(C, P) 
    \end{equation}
    is isomorphism. On the right we have isomorphisms
    \begin{eqnarray*}
        (\text{RHS}) &\simeq & \Cone\bigl( \ul{\Hom}_{\sO_{X,\bbox}}(f_!f^!M\oplus i_!i^!M,P) \to   \ul{\Hom}_{\sO_{X,\bbox}}(f_!j_!j^!f^! M,P)\bigr)[-1]\\
        &\overset{\text{Thm}~\ref{CSex}~(1), (7)}{\simeq} & \Cone\bigl( \ul{\Hom}_{\sO_{X,\bbox}}(M,f_*f^*P\oplus i_*i^*P) \to   \ul{\Hom}_{\sO_{X,\bbox}}(M,f_*j_*j^*f^*P)\bigr)[-1]\\
        &\simeq &\ul{\Hom}_{\sO_{X,\bbox}}\bigl(M,\Cone( f_*f^*P\oplus i_*i^*P \to f_*j_*j^*f^*P)\bigr)[-1].
    \end{eqnarray*}
    Here, by the projection formula Theorem~\ref{CSex} (3) and Nisnevich descent for coherent sheaves, we have isomorphisms
    \[
    \Cone( f_*f^*P\oplus i_*i^*P \to f_*j_*j^*f^*P) \simeq P\otimes_{\sO_{X,\bbox}}^{\L} \Cone( f_*f^*\sO_X\oplus i_*i^*\sO_X \to f_*j_*j^*f^*\sO_X) \simeq P[1].
    \]
    Thus, we obtain the claim.
\end{proof}

 Theorem~\ref{main} (8) follows from Theorem~\ref{main} (1), (2), (3) and (7).

\end{Observation}

\section{\texorpdfstring{$m$}{m}-Solid realization of modulus pairs}

Let $(X,D)$ be a normal modulus pair over $A$. In the previous section, we defined $\Solid (X,D)$. In this section, we define for $m\geq 0$:
\[
\Solidm(X,D)
\]
in $D(A_{\bbox})$ when $X$ is regular, and $D$ is a strictly normal crossing effective Cartier divisor. The construction is characterized by the property that its dual $ R\ul{\Hom}_{A_\bbox}(\Solidm(X, D), A)$ is discrete, and that the $i$-th cohomology group of this dual coincides with Kelly–Miyazaki’s Hodge cohomology with modulus:
\[
H^i\bigl((\Solidm(X,D))^\vee\bigr)\;\simeq\; H_\Zar^i(X,\Omega^m(\log |D|)(D-|D|)).
\]

If $X$ is regular and $|D|$ is strictly normal crossing, then we call a pair $(X,D)$ \textit{regular modulus pair}. For a regular modulus pair $(X,D)$ over $A$, we define:
\[
\Solidm(X,D)=f_! \ul{\Hom}_{\sO_{X,\bbox}} \Bigl(\Omega_{X/A}^m(\log |D|), f^!A\otimes_{\sO_{X,\bbox}} \sO_X(-D+|D|) \Bigr).
\]
For $m=0$, $\Solid^0(X,D)=\Solid(X,D)$. When $X$ is smooth over $A$, then we have
\begin{eqnarray*}
\Solidm(X,D) &=& f_! \Omega^{\dim X/A-m}_{X/A}(\log |D|)(-D)[\dim X/A].
\end{eqnarray*}

$\Solidm(-)$ satisfies the following properties similar to $\Solid$:
\begin{thm}\label{mSolidthm}
     Let $A$ be a discrete Noetherian commutative ring. Let $A_{\bbox}$ be the analytic ring $(A, A)_{\bbox}$, and let $D(A_{\bbox})$ be the derived category of solid $A_{\bbox}$-modules. 
    \begin{itemize}
        \item[(1)] Let $h: (Y, C) \to (X, D)$ be an ambient morphism of regular modulus pairs over $A$. Then, there exists a morphism in $D(A_{\bbox})$:
        \[
        h_+^m: \Solidm(Y, C) \to \Solidm(X, D).
        \]
        Furthermore, for a sequence of ambient morphisms $(Z, B) \overset{k}{\to} (Y, C) \overset{h}{\to} (X, D)$, the functoriality $(h \circ k)_+^m = h_+^m \circ k_+^m$ holds.
        \item[(2)] $\bcube$-invariance holds. That is, the natural morphism:
        \[
        \Solidm((X, D) \otimes (\PP^1,\{\infty\})) \to \Solidm(X, D)
        \]
        becomes an isomorphism in $D(A_{\bbox})$.
        \item[(3)]  Assume $X \to \Spec A$ is locally complete intersection. Given an integral regular scheme $Z$ and a regular closed embedding $Z \hookrightarrow X$ with $|Z| \subset |D|$ and assume $|p^*(D)|$ is SNCD, then the invariance under blow-ups holds. That is, the natural morphism:
        \[
        \Solidm(\text{Bl}_Z X, p^* D) \to \Solidm(X, D)
        \]
        becomes an isomorphism. Here, $p: \text{Bl}_Z X \to X$ denotes the blow-up of $X$ along $Z$.
        \item[(4)] Assume that $X$ is proper smooth over $A$. Then, $\Solidm(X, D)$ is a discrete solid $A_{\bbox}$-module complex (that is a usual $A$- module complex), and its $i$-th cohomology group is:
        \[
        H^{i + \dim X/A}(X, \Omega_{X/A}^{\dim X/A-m}(\log |D|)(-D)).
        \]
        When $A$ is a field, by Serre duality, this coincides with:
        \[
        H^{-i}(X, \Omega^m_{X/A}(\log |D|)(D-|D|))^{\vee}.
        \]
        \item[(5)] Assume $X$ is smooth over $A$. The dual of $\Solidm(X, D)$ over $A$:
        \[
        R\ul{\Hom}_{A_\bbox}(\Solidm(X, D), A)
        \]
        is a discrete solid $A_{\bbox}$-module, and it is quasi-isomorphic to 
        \[
        R\Gamma_{\Zar}(X,\Omega^m(\log |D|)(D-|D|)).
        \]
        \item[(6)] For a $\underline{MV}$-square
        \[
        \xymatrix{
        (V,j^*C) \ar[r]^-{j} \ar[d]_{g} & (Y,C) \ar[d]^{f}\\
        (U,i^*D) \ar[r]^-{i} & (X,D)
        }
        \]
        there exists a distinguished triangle 
        \[
        \Solidm(V,j^*C) \overset{(j^m_+,g^m_+)}{\to} \Solidm(Y,C) \oplus \Solidm(U,i^*D) \overset{i^m_+-f^m_+}{\to} \Solidm(X,D) \overset{+}{\to} 
        \]

        \item[(7)] We assume $A$ is a field of characteristic $0$. There is a covariant functor 
        \[
        \ulMDMgm(A,\Q) \to D(A_{\bbox}) 
        \]
        which sends $\ulM(X,D)$ to $\Solidm(X,D)$.
        \end{itemize}
        \end{thm} We begin by proving the following lemma.
            \begin{lemma}\label{compactprojection}
                Let $h:Y \to X$ be a locally of finite type morphism between separated Noetherian schemes. Let $P$ be a bounded perfect discrete complex on $X$. For any $Q\in D(\sO_{X,\bbox})$, there is an isomorphism
                \[
                h^!\ul{\Hom}_{\sO_{X,\bbox}}(P,Q)\simeq \ul{\Hom}_{\sO_{Y,\bbox}}(h^*P,h^!Q)
                \]
                in $D(\sO_{Y,\bbox})$.\end{lemma}
                    \begin{proof} We have isomorphisms 
                    \begin{eqnarray*}
        h^!\ul{\Hom}_{\sO_{X,\bbox}}(P,Q) \simeq h^!\sO_{X}\otimes^\L_{\sO_{Y,\bbox}} h^*\ul{\Hom}_{\sO_{X,\bbox}}(P,Q)\\
        \simeq h^!\sO_{X}\otimes^\L_{\sO_{Y,\bbox}} \ul{\Hom}_{\sO_{Y,\bbox}}(h^*P,h^*Q)
                    \end{eqnarray*}
                    Since $P$ is dualizable in $D(\sO_{X,\bbox})$, we obtain 
         \begin{eqnarray*}
             h^!\sO_{X}\otimes^\L_{\sO_{Y,\bbox}} \ul{\Hom}_{\sO_{Y,\bbox}}(h^*P,h^*Q) &\simeq& h^!\sO_{X}\otimes^\L_{\sO_{Y,\bbox}}(h^*P)^{\vee}\otimes^\L_{\sO_{Y,\bbox}}h^*Q\\
             &\simeq & \ul{\Hom}_{\sO_{Y,\bbox}}(h^*P,h^!\sO_X \otimes^\L_{\sO_{Y,\bbox}} h^*Q) \\
             &\simeq &\ul{\Hom}_{\sO_{Y,\bbox}}(h^*P,h^!Q).
         \end{eqnarray*}
         This is what we wanted to prove.
        \end{proof}
Let $A$ be a discrete Noetherian commutative ring. Let $h: (Y, C) \to (X, D)$ be an ambient morphism of regular modulus pairs over $A$. Write $f:X \to \Spec A$ and $g:Y \to \Spec A$ for the structure morphisms. By \cite{kelly2023hodge}, the map $h$ induces a morphism of Hodge sheaves with modulus $h^*\Omega_{X/A}^m(\log |D|)(D-|D|) \to \Omega_{Y/A}^m(\log |C|)(C-|C|)$.
By contravariance of the internal Hom in the first variable, this yields a morphism
\[
\ul{\Hom}_{\sO_{Y,\bbox}} \Bigl(\Omega_{Y/A}^m(\log |C|)(C-|C|), g^!A \Bigr) \to  \ul{\Hom}_{\sO_{Y,\bbox}}\Bigl(h^*\Omega_{X/A}^m(\log |D|)(D-|D|), g^!A \Bigr).
\]
Moreover, by Lemma~\ref{compactprojection} there is a canonical morphism
\begin{equation}\label{solidmmorphism1}
\ul{\Hom}_{\sO_{Y,\bbox}} \Bigl(\Omega_{Y/A}^m(\log |C|)(C-|C|), g^!A \Bigr) \to  h^!\ul{\Hom}_{\sO_{X,\bbox}}\Bigl(\Omega_{X/A}^m(\log |D|)(D-|D|), f^!A \Bigr).
\end{equation}
Applying $h_!$ and composing with the counit $h_!h^! \Longrightarrow id$, we obtain a morphism
\[
h^m_+:\Solidm(Y,C) \to \Solidm(X,D).
\]
By the functoriality of Hodge sheaves with modulus (see \cite{kelly2023hodgecohomologyramificationfiltration}), this construction satisfies functoriality. In the case $m=0$, the above construction agrees with that of Section~2.

\begin{proof}[Proof of Theorem~\ref{mSolidthm} (2) (resp. (3))]
    Let $f:X\to\Spec A$ and $g:Y\to\Spec A$ denote the structure morphisms, where
\[
Y:=X\times_A\PP^1 \ \ (\text{resp. } Y:=\Bl_Z X),\qquad
C:=X\times_A\{\infty\}+D\times_A\PP^1 \ \ (\text{resp. } C:=p^{*}D)
\]
is an effective Cartier divisor on $Y$. We first show that the morphism
\begin{equation}\label{mSolidmorphism2}
p_!\,\ul{\Hom}_{\sO_{Y,\bbox}}\!\Bigl(\Omega_{Y/A}^m(\log|C|)(C-|C|),\,g^!A\Bigr)
\xrightarrow{\,p_!\eqref{solidmmorphism1}\,}
p_!\,\ul{\Hom}_{\sO_{Y,\bbox}}\!\Bigl(p^*\Omega_{X/A}^m(\log|D|)(D-|D|),\,g^!A\Bigr)
\end{equation}
is an isomorphism.
In $D(\Coh(Y))$, let $G$ denote the cone of the morphism
$\,p^*\Omega_{X/A}^m(\log|D|)(D-|D|)\to\Omega_{Y/A}^m(\log|C|)(C-|C|)$.
By \cite[Cor.~5.2]{kelly2023hodgecohomologyramificationfiltration} (resp. \cite[Thm.~6.1]{kelly2023hodgecohomologyramificationfiltration}), we have $Rp_*G=0$.
Hence, by Verdier duality,
\[
p_!\,\ul{\Hom}_{\sO_{Y,\bbox}}\!\bigl(G,\,g^!A\bigr)\;\simeq\;
\ul{\Hom}_{\sO_{X,\bbox}}\!\bigl(p_*G,\,f^!A\bigr)\;=\;0,
\]
and \eqref{mSolidmorphism2} is an isomorphism.
Therefore, to establish the claim, it suffices to show that the counit is an isomorphism:
\[
p_!p^!\ul{\Hom}_{\sO_{X,\bbox}}\Bigl(\Omega_{X/A}^m(\log |D|)(D-|D|), f^!A \Bigr) \to \ul{\Hom}_{\sO_{X,\bbox}}\Bigl(\Omega_{X/A}^m(\log |D|)(D-|D|), f^!A \Bigr).
\]
Set
\[
F:=\Omega_{X/A}^m(\log|D|)(D-|D|).
\]
By Lemma~\ref{compactprojection}, Verdier duality, and the projection formula, the left-hand side identifies as
\begin{eqnarray*}
    p_!p^!\,\ul{\Hom}_{\sO_{X,\bbox}}(F, f^!A)
&\overset{\text{Lemma~\ref{compactprojection}}}{\simeq}&
p_!\,\ul{\Hom}_{\sO_{Y,\bbox}}(p^*F, g^!A) \\
&\overset{\text{Verdier duality}}{\simeq}&
\ul{\Hom}_{\sO_{X,\bbox}}(Rp_*p^*F, f^!A)
\\
&\overset{\text{projection formula}}{\simeq}&
\ul{\Hom}_{\sO_{X,\bbox}}\bigl(F\otimes^\L_{\sO_{X,\bbox}} Rp_*\sO_{Y},\, f^!A\bigr)\\
&\overset{Rp_*\sO_{Y}\simeq \sO_{X}.}{\simeq}&
\ul{\Hom}_{\sO_{X,\bbox}}(F, f^!A)
\end{eqnarray*}
Consequently, the counit map above is an isomorphism.
\end{proof}

Theorem~\ref{mSolidthm} (4) can be directly proven using Theorem~\ref{CSex} (4), (5).

\begin{Observation}[Proof of Theorem~\ref{mSolidthm} (5)]
Assume $f:X\to \Spec A $ is smooth. By Theorem~\ref{CSex} (1),
\[
R\ul{\Hom}_{A_\bbox}(\Solidm(X, D), A) = f_* \ul{\Hom}_{\sO_{X,\bbox}}\Bigl(\ul{\Hom}_{\sO_{X,\bbox}} \bigl(\Omega_{X/A}^m(\log |D|)(D-|D|), f^!A\bigr),f^!A\Bigr).
\]
By Theorem~\ref{CSex} (6), $f^!A$ is invertible and $\Omega_{X/A}^m(\log |D|)(D-|D|)$ is dualizable, we obtain the claim
\[
R\ul{\Hom}_{A_\bbox}(\Solidm(X, D), A) = f_* \Omega_{X/A}^m(\log |D|)(D-|D|).
\]
\end{Observation}

Theorem~\ref{mSolidthm} (6) can be directly proven using Theorem~\ref{CSex} (7) and proposition~\ref{solidNisdescent}.

 Theorem~\ref{mSolidthm} (7) follows from Theorem~\ref{mSolidthm} (1), (2), (3) and (6).

\section{Explicit formulas using models of finite type over \texorpdfstring{$\Z$}{Z}}
\label{sec:boundary}

In this section, we give an explicit description of $\Solidm(X,D)$ for a smooth variety $X$ over a field $A$ (Theorem~\ref{prosolidthm}). As an application, for an open immersion $U\hookrightarrow Y$, we compute the image of the cone of $\ulM(U,D|_U)\to\ulM(Y,D)$ under $\Solidm$ (Theorem~\ref{solidboundarycone}).

The main tool is Proposition~\ref{renprosolid} below. It describes the functor $j_!j^*$ for an open immersion $j$ of schemes of finite type over $\Z$ in terms of the formal completion along the closed complement. For a principal open subset defined by a non-zero-divisor, this is Ren's comparison theorem \cite[Theorem~3.3]{prosolid}. In general, the field $A$ is not of finite type over $\Z$. Therefore, we take the completion on a model over a subring $A_0\subset A$ of finite type over $\Z$, and then take the base change to $A$ (Proposition~\ref{solidboundary}). We note that the completion computed directly over $A$ gives a different object in general (Remark~\ref{rem:naive}).

Throughout this section, for a discrete commutative ring $B$, we write $B_{\bbox}=(B,B)_{\bbox}$. All functors are derived functors between derived categories of solid modules. For a ring homomorphism $A_0\to A$, the base change functor $D(A_{0,\bbox})\to D(A_{\bbox})$ is given by $M\mapsto M\otimes^{\L}_{A_{0,\bbox}}A_{\bbox}$. For schemes over a ring which is not of finite type over $\Z$, we use the six-functor formalism of \cite[Theorem~9.15]{ScholzeSix} (see also Remark~\ref{remark} (2)).

Let $Y_0$ be a separated scheme of finite type over $\Z$, and let $Z_0\subset Y_0$ be a closed subscheme defined by an ideal sheaf $\sI_0$. For a perfect complex $F_0$ on $Y_0$, we define the \textit{completion of $F_0$ along $Z_0$} by
\begin{equation}\label{model:completion}
\widehat F_0:=R\varprojlim_{n\geq1}\bigl(F_0\otimes^{\L}_{\sO_{Y_0}}\sO_{Y_0}/\sI_0^n\bigr)\in D(\sO_{Y_0,\bbox}).
\end{equation}
Here, each $F_0\otimes^{\L}_{\sO_{Y_0}}\sO_{Y_0}/\sI_0^n$ is regarded as a discrete solid $\sO_{Y_0,\bbox}$-module, and the limit is taken in $D(\sO_{Y_0,\bbox})$. The closed subscheme $Z_0$ will always be clear from the context. If $Y_0=\Spec B$ is affine, $I=\Gamma(Y_0,\sI_0)$, and $F_0$ is represented by a bounded complex $P$ of finitely generated projective $B$-modules, then $\widehat F_0$ is represented by the complex $\varprojlim_nP/I^nP$ of condensed $B$-modules, since the transition maps are surjective. In particular, if $F_0$ is a finite locally free sheaf, then $\widehat F_0$ is the usual formal completion of $F_0$ along $Z_0$ with its $\sI_0$-adic topology.

\begin{prop}\label{renprosolid}
Let $Y_0$, $Z_0$ and $F_0$ be as above, and let $j_0:U_0:=Y_0\setminus Z_0\hookrightarrow Y_0$ be the open immersion. Then there exists a distinguished triangle in $D(\sO_{Y_0,\bbox})$:
\begin{equation}\label{model:localization}
j_{0!}j_0^*F_0\longrightarrow F_0\longrightarrow\widehat F_0\overset{+}{\longrightarrow}
\end{equation}
where the first morphism is the counit of the adjunction $j_{0!}\dashv j_0^*$, and the second morphism is the natural map. In particular, $\widehat F_0$ depends only on the open subscheme $U_0$, and not on the ideal sheaf $\sI_0$.
\end{prop}

\begin{proof}
First, we construct a morphism from $\Cone(j_{0!}j_0^*F_0\to F_0)$ to $\widehat F_0$. For every $n$, we have $j_0^*(F_0\otimes^{\L}_{\sO_{Y_0}}\sO_{Y_0}/\sI_0^n)=0$. Since $j_0^*$ has the left adjoint $j_{0!}$, it commutes with limits, and hence $j_0^*\widehat F_0=0$. Thus, for every $k\in\Z$, we have
\[
\Hom\bigl(j_{0!}j_0^*F_0,\widehat F_0[k]\bigr)\simeq\Hom\bigl(j_0^*F_0,j_0^*\widehat F_0[k]\bigr)=0,
\]
and the natural map $F_0\to\widehat F_0$ factors uniquely through a morphism
\[
c:\Cone(j_{0!}j_0^*F_0\to F_0)\longrightarrow\widehat F_0.
\]
It suffices to show that $c$ is an isomorphism. Let $k:V\hookrightarrow Y_0$ be an affine open subscheme. The functor $k^*$ commutes with $j_{0!}$ by Theorem~\ref{CSex} (7), with limits since it has the left adjoint $k_!$, and with tensor products. Hence $k^*(c)$ is the morphism $c$ for $(V,Z_0\cap V,k^*F_0)$. By solid Zariski descent \cite[Corollary~5.4.17]{CamargoSolid}, we may assume that $Y_0=\Spec B$ is affine and that $\sI_0$ corresponds to an ideal $I=(t_1,\ldots,t_r)\subset B$.

For $t\in B$, let $j_t:\Spec B[1/t]\hookrightarrow\Spec B$ be the open immersion. By \cite[Lemmas~5.4.13 (2), 5.4.14 and the proof of Proposition~5.4.15]{CamargoSolid}, there is a distinguished triangle in $D(B_{\bbox})$:
\begin{equation}\label{model:principal}
j_{t!}j_t^*B\longrightarrow B\longrightarrow\widehat B_{(t)}\overset{+}{\longrightarrow}
\end{equation}
where $\widehat B_{(t)}=\varprojlim_nB/t^nB$ with its $t$-adic topology, the first morphism is the counit, and the second morphism is the natural map. If $t$ is a non-zero-divisor, \eqref{model:principal} is also obtained from \cite[Theorem~3.3 and (3.3.2)]{prosolid}. Let $M\in D(B_{\bbox})$. By the projection formula (Theorem~\ref{CSex} (3)), tensoring \eqref{model:principal} with $M$ gives a distinguished triangle
\[
j_{t!}j_t^*M\longrightarrow M\longrightarrow M\otimes^{\L}_{B_{\bbox}}\widehat B_{(t)}\overset{+}{\longrightarrow}.
\]
Since $j_{t!}$ is fully faithful, we have $j_t^*j_{t!}\simeq\id$, and hence $j_t^*(M\otimes^{\L}_{B_{\bbox}}\widehat B_{(t)})=0$.

Now, let $j:U:=\Spec B\setminus V(I)\hookrightarrow\Spec B$ be the open immersion. By \cite[Proposition~5.3.12 (2)]{CamargoSolid}, we have
\[
\widehat B_{(t_1)}\otimes^{\L}_{B_{\bbox}}\cdots\otimes^{\L}_{B_{\bbox}}\widehat B_{(t_r)}\simeq\widehat B_I:=\varprojlim_nB/I^n.
\]
We show that for any $M\in D(B_{\bbox})$, there is a distinguished triangle
\begin{equation}\label{model:affine}
j_!j^*M\longrightarrow M\longrightarrow M\otimes^{\L}_{B_{\bbox}}\widehat B_I\overset{+}{\longrightarrow}
\end{equation}
whose first morphism is the counit. Set $N:=M\otimes^{\L}_{B_{\bbox}}\widehat B_I$. By the above, we have $j_{t_i}^*N=0$ for each $i$. Since $U=\bigcup_i\Spec B[1/t_i]$, Zariski descent gives $j^*N=0$. On the other hand, the morphism $M\to N$ is the composite
\[
M\to M\otimes^{\L}_{B_{\bbox}}\widehat B_{(t_1)}\to M\otimes^{\L}_{B_{\bbox}}\widehat B_{(t_1)}\otimes^{\L}_{B_{\bbox}}\widehat B_{(t_2)}\to\cdots\to N,
\]
and by the above, the fiber of each morphism is of the form $j_{t_i!}(-)$. Since $j_{t_i}$ factors through $j$, these fibers lie in the essential image of $j_!$. The functor $j_!$ is exact and fully faithful, so its essential image is closed under extensions. Hence $K:=\fib(M\to N)$ lies in the essential image of $j_!$, and the counit $j_!j^*K\to K$ is an isomorphism. Since $j^*N=0$, we have $j^*K\simeq j^*M$. Therefore, $K\simeq j_!j^*M$, and by the naturality of the counit, the morphism $K\to M$ is the counit. This proves \eqref{model:affine}.

Finally, let $M=F_0$. Since $F_0$ is perfect, the functor $F_0\otimes^{\L}_{B_{\bbox}}(-)$ commutes with limits. Hence we have
\[
F_0\otimes^{\L}_{B_{\bbox}}\widehat B_I\simeq R\varprojlim_n\bigl(F_0\otimes^{\L}_BB/I^n\bigr)=\widehat F_0,
\]
and the second morphism in \eqref{model:affine} is the natural map $F_0\to\widehat F_0$. Therefore, $c$ is an isomorphism. The last assertion follows from \eqref{model:localization}, since the counit $j_{0!}j_0^*F_0\to F_0$ depends only on $U_0$.
\end{proof}

\begin{prop}\label{solidboundary}
Let $A$ be a field, and let $A_0\subset A$ be a subring of finite type over $\Z$. Let $g_0:Y_0\to\Spec A_0$ be a proper morphism, let $Z_0\subset Y_0$ be a closed subscheme with ideal sheaf $\sI_0$, and let $F_0$ be a perfect complex on $Y_0$. Consider the following Cartesian diagram:
\[
\xymatrix{
U \ar[r]^-{j} \ar[d] & Y \ar[r]^-{g} \ar[d]^{q} & \Spec A \ar[d]\\
U_0 \ar[r]_-{j_0} & Y_0 \ar[r]_-{g_0} & \Spec A_0
}
\]
where $U_0=Y_0\setminus Z_0$. Set $F=q^*F_0$, $f=g\circ j$ and $f_0=g_0\circ j_0$. Then there exists a distinguished triangle in $D(A_{\bbox})$:
\begin{equation}\label{model:pushforward-triangle}
f_!j^*F\longrightarrow g_*F\longrightarrow\bigl(g_{0*}\widehat F_0\bigr)\otimes^{\L}_{A_{0,\bbox}}A_{\bbox}\overset{+}{\longrightarrow}
\end{equation}
where the first morphism is induced by the counit $j_!j^*F\to F$, and $\widehat F_0$ is the completion of $F_0$ along $Z_0$. In particular, we have
\begin{equation}\label{model:pushforward-fiber}
f_!j^*F\simeq\fib\Bigl(g_*F\longrightarrow\bigl(g_{0*}\widehat F_0\bigr)\otimes^{\L}_{A_{0,\bbox}}A_{\bbox}\Bigr).
\end{equation}
The object $(g_{0*}\widehat F_0)\otimes^{\L}_{A_{0,\bbox}}A_{\bbox}$, together with its morphism from $g_*F$, does not depend on the choice of $(A_0,Y_0,Z_0,F_0)$.

Moreover, let $i_n:Z_{0,n}:=V(\sI_0^n)\hookrightarrow Y_0$ be the closed immersion and set $g_n:=g_0\circ i_n$. Then we have
\begin{equation}\label{model:formal-pushforward}
g_{0*}\widehat F_0\simeq R\varprojlim_{n\geq1}Rg_{n*}(Li_n^*F_0).
\end{equation}
Here, each $Rg_{n*}(Li_n^*F_0)$ is the usual derived pushforward. It is a bounded complex with finitely generated cohomology groups, and we regard it as a discrete solid $A_{0,\bbox}$-module.
\end{prop}

\begin{proof}
Since $g_0$ is proper, we have $g_{0!}=g_{0*}$ by Theorem~\ref{CSex} (5). Applying $g_{0*}$ to \eqref{model:localization}, we obtain a distinguished triangle in $D(A_{0,\bbox})$:
\begin{equation}\label{model:triangle-A0}
f_{0!}j_0^*F_0\longrightarrow g_{0*}F_0\longrightarrow g_{0*}\widehat F_0\overset{+}{\longrightarrow}.
\end{equation}
Since $A_0$ is a subring of the field $A$, $A_0$ is an integral domain, and $A$ is flat over $A_0$ (it is flat over the fraction field of $A_0$). Hence the squares in the above diagram are also derived Cartesian. By base change \cite[Theorem~9.15]{ScholzeSix} (cf.\ Theorem~\ref{CSex} (7)), we have isomorphisms
\begin{eqnarray*}
(g_{0*}F_0)\otimes^{\L}_{A_{0,\bbox}}A_{\bbox}&\simeq&g_*q^*F_0=g_*F,\\
(f_{0!}j_0^*F_0)\otimes^{\L}_{A_{0,\bbox}}A_{\bbox}&\simeq&f_!j^*q^*F_0=f_!j^*F,
\end{eqnarray*}
which are compatible with the counits. Applying $-\otimes^{\L}_{A_{0,\bbox}}A_{\bbox}$ to \eqref{model:triangle-A0}, we obtain \eqref{model:pushforward-triangle} and \eqref{model:pushforward-fiber}. By \eqref{model:pushforward-triangle}, the object $(g_{0*}\widehat F_0)\otimes^{\L}_{A_{0,\bbox}}A_{\bbox}$ is the cone of $f_!j^*F\to g_*F$, which is defined without the model. This proves the independence of the choice.

Next, we prove \eqref{model:formal-pushforward}. By the projection formula, we have $F_0\otimes^{\L}_{\sO_{Y_0}}\sO_{Y_0}/\sI_0^n\simeq i_{n*}Li_n^*F_0$. Since $g_{0*}$ is a right adjoint, it commutes with limits. Hence we obtain
\[
g_{0*}\widehat F_0\simeq R\varprojlim_n g_{0*}i_{n*}Li_n^*F_0=R\varprojlim_n g_{n*}Li_n^*F_0.
\]
Since $g_n$ is proper and $Li_n^*F_0$ is a discrete bounded complex, we have $g_{n*}(Li_n^*F_0)=Rg_{n*}(Li_n^*F_0)$ by Theorem~\ref{CSex} (5). This is what we wanted to prove.
\end{proof}

\begin{remark}\label{rem:naive}
In \eqref{model:pushforward-triangle}, the base change $-\otimes^{\L}_{A_{0,\bbox}}A_{\bbox}$ is taken after the limit in \eqref{model:formal-pushforward}. These two operations do not commute in general. Hence the third term of \eqref{model:pushforward-triangle} is in general different from $R\varprojlim_nR\Gamma(Z_n,F|_{Z_n})$, where $Z_n\subset Y$ is the closed subscheme defined by $\sI_0^n\sO_Y$. For example, let $A=\Q$, $A_0=\Z$, $Y_0=\PP^1_{\Z}$, let $Z_0=\{x=0\}$ be the zero section, and let $F_0=\sO_{Y_0}$. Then $g_{0*}\widehat F_0=\Z[[x]]\simeq\prod_{\N}\Z$. By the description of the free $\Q_{\bbox}$-modules \cite[Warning~9.11]{ScholzeSix}, we have
\[
\Z[[x]]\otimes^{\L}_{\Z_{\bbox}}\Q_{\bbox}\simeq\bigcup_{N\geq1}\Z[1/N][[x]],
\]
that is, the ring of power series with bounded denominators. This is strictly smaller than $\Q[[x]]=\varprojlim_n\Q[x]/(x^n)$. For example, $\sum_{n}x^n/n!$ is not contained in it.
\end{remark}

We now apply Proposition~\ref{solidboundary} to $\Solidm(X,D)$. We set $\Omega^r=0$ for $r<0$.

\begin{thm}\label{prosolidthm}
Let $A$ be a field, let $X$ be a smooth separated scheme of finite type over $A$ of pure dimension $d$, and let $D$ be an effective Cartier divisor on $X$. Assume that there exist an open immersion $j:X\hookrightarrow Y$ into a smooth proper $A$-scheme $Y$ of pure dimension $d$ and an effective Cartier divisor $\widetilde D$ on $Y$ such that $\widetilde D|_X=D$ and $|\widetilde D|$ is a simple normal crossings divisor. Let $g:Y\to\Spec A$ be the structure morphism, let $m\geq0$, and set
\[
Z:=Y\setminus X,\qquad F_m:=\Omega^{d-m}_{Y/A}(\log|\widetilde D|)(-\widetilde D).
\]
Then the following assertions hold.
\begin{itemize}
\item[(1)] There exist a subring $A_0\subset A$ of finite type over $\Z$, a proper morphism $g_0:Y_0\to\Spec A_0$ with $Y_0\times_{A_0}A\simeq Y$, a closed subscheme $Z_0\subset Y_0$ with $Z_0\times_{A_0}A=Z$, and a finite locally free $\sO_{Y_0}$-module $F_{m,0}$ whose pullback to $Y$ is isomorphic to $F_m$.
\item[(2)] For any such $(A_0,Y_0,Z_0,F_{m,0})$, there is an isomorphism in $D(A_{\bbox})$:
\begin{equation}\label{model:realization-fiber}
\Solidm(X,D)\simeq\fib\Bigl(g_*F_m\longrightarrow\bigl(g_{0*}\widehat F_{m,0}\bigr)\otimes^{\L}_{A_{0,\bbox}}A_{\bbox}\Bigr)[d],
\end{equation}
where $\widehat F_{m,0}$ is the completion of $F_{m,0}$ along $Z_0$. The right-hand side does not depend on the choice of $(A_0,Y_0,Z_0,F_{m,0})$.
\end{itemize}
\end{thm}

If $A$ has characteristic $0$ and $|D|$ is a simple normal crossings divisor on $X$, such $Y$ and $\widetilde D$ exist by Nagata's compactification theorem and resolution of singularities.

\begin{proof}
First, we prove (1). Since $Y$ is proper over $A$ and $A$ is the filtered union of its subrings of finite type over $\Z$, there exist $A_0$ and a proper $A_0$-scheme $Y_0$ with $Y_0\times_{A_0}A\simeq Y$ \cite[Tags~01ZM and 081F]{Stacks}. Let $q:Y\to Y_0$ be the projection, and let $\sI\subset\sO_Y$ be the ideal sheaf of $Z$. After enlarging $A_0$, the inclusion $\sI\hookrightarrow\sO_Y$ is the pullback of a morphism $\sI_0'\to\sO_{Y_0}$ of $\sO_{Y_0}$-modules of finite presentation \cite[Tag~01ZR]{Stacks}, and $F_m$ is the pullback of a finite locally free $\sO_{Y_0}$-module $F_{m,0}$ \cite[Tag~0B8W]{Stacks}. Let $Z_0$ be the closed subscheme of $Y_0$ defined by the image of $\sI_0'\to\sO_{Y_0}$. Since $A$ is flat over $A_0$, the functor $q^*$ preserves images, and hence $Z_0\times_{A_0}A=Z$.

Next, we prove (2). Let $f=g\circ j$. By Theorem~\ref{CSex} (4), we have $f^!A\simeq\Omega^d_{X/A}[d]$. Since $|D|$ is a simple normal crossings divisor, $\Omega^1_{X/A}(\log|D|)$ is locally free of rank $d$, and the wedge product gives a perfect pairing
\[
\Omega^m_{X/A}(\log|D|)\otimes_{\sO_X}\Omega^{d-m}_{X/A}(\log|D|)\longrightarrow\Omega^d_{X/A}(\log|D|)=\Omega^d_{X/A}(|D|).
\]
Hence, we obtain the following isomorphisms:
\begin{eqnarray*}
\Solidm(X,D)&=&f_!\ul{\Hom}_{\sO_{X,\bbox}}\Bigl(\Omega^m_{X/A}(\log|D|),f^!A\otimes^{\L}_{\sO_{X,\bbox}}\sO_X(-D+|D|)\Bigr)\\
&\overset{\text{Thm~\ref{CSex} (4)}}{\simeq}&f_!\Bigl(\Omega^m_{X/A}(\log|D|)^{\vee}\otimes_{\sO_X}\Omega^d_{X/A}(-D+|D|)\Bigr)[d]\\
&\overset{\text{wedge product}}{\simeq}&f_!\,\Omega^{d-m}_{X/A}(\log|D|)(-D)[d]\\
&=&f_!j^*F_m[d].
\end{eqnarray*}
Applying Proposition~\ref{solidboundary} to $(Y_0,Z_0,F_{m,0})$, we obtain \eqref{model:realization-fiber}. The independence of the choice also follows from Proposition~\ref{solidboundary}.
\end{proof}

\begin{remark}\label{rem:logmodel}
In Theorem~\ref{prosolidthm}, we can choose $F_{m,0}$ as a sheaf of logarithmic differential forms on a model. Write $\widetilde D=\sum_{\alpha}n_{\alpha}H_{\alpha}$, where $H_{\alpha}$ are the irreducible components of $|\widetilde D|$. After enlarging $A_0$, we may assume that $Y_0$ is smooth over $A_0$ \cite[Tag~0C0C]{Stacks}, that each $H_{\alpha}$ is the base change of a closed subscheme $H_{\alpha,0}\subset Y_0$, and that for every set $J$ of indices, $\bigcap_{\alpha\in J}H_{\alpha,0}$ is either empty or smooth over $A_0$ of pure relative dimension $d-|J|$. Then $H_0:=\sum_{\alpha}H_{\alpha,0}$ is a relative simple normal crossings divisor on $Y_0$ over $A_0$. Hence $\Omega^r_{Y_0/A_0}(\log H_0)$ is finite locally free, its formation commutes with base change, and we can take
\[
F_{m,0}=\Omega^{d-m}_{Y_0/A_0}(\log H_0)(-\widetilde D_0),\qquad\widetilde D_0:=\sum_{\alpha}n_{\alpha}H_{\alpha,0}.
\]
For $m=0$, we have $\Solid^0=\Solid$ (see Section~3), and $F_{0,0}=\omega_{Y_0/A_0}(-\widetilde D_0+H_0)$. Thus, Theorem~\ref{prosolidthm} also describes the solid realization $\Solid(X,D)$ of Section~2.
\end{remark}

\begin{thm}[The cone of an open immersion]\label{solidboundarycone}
Let $A$ be a field of characteristic $0$, let $Y$ be a smooth proper $A$-scheme of pure dimension $d$, and let $D$ be an effective Cartier divisor on $Y$ such that $|D|$ is a simple normal crossings divisor. Let $j:U\hookrightarrow Y$ be an open immersion, let $m\geq0$, and set
\[
Z:=Y\setminus U,\qquad F_m:=\Omega^{d-m}_{Y/A}(\log|D|)(-D).
\]
Choose a model $(A_0,Y_0,Z_0,F_{m,0})$ of $(Y,Z,F_m)$ as in Theorem~\ref{prosolidthm} (1), and let $g_0:Y_0\to\Spec A_0$ be the structure morphism. Then there is an isomorphism in $D(A_{\bbox})$:
\begin{equation}\label{model:motivic-cone}
\Solidm\Bigl(\Cone\bigl(\ulM(U,D|_U)_{\Q}\longrightarrow\ulM(Y,D)_{\Q}\bigr)\Bigr)\simeq\Bigl(\bigl(g_{0*}\widehat F_{m,0}\bigr)\otimes^{\L}_{A_{0,\bbox}}A_{\bbox}\Bigr)[d].
\end{equation}
In particular, the right-hand side does not depend on the choice of $(A_0,Y_0,Z_0,F_{m,0})$.
\end{thm}

\begin{proof}
Let $g:Y\to\Spec A$ be the structure morphism, and let $f=g\circ j$. Set
\[
G:=\ul{\Hom}_{\sO_{Y,\bbox}}\Bigl(\Omega^m_{Y/A}(\log|D|),g^!A\otimes^{\L}_{\sO_{Y,\bbox}}\sO_Y(-D+|D|)\Bigr),
\]
so that $\Solidm(Y,D)=g_!G$. As in the proof of Theorem~\ref{prosolidthm}, we have $G\simeq F_m[d]$. Since $j^!=j^*$ and $j^*\Omega^m_{Y/A}(\log|D|)=\Omega^m_{U/A}(\log|D|_U|)$, Lemma~\ref{compactprojection} gives
\[
\Solidm(U,D|_U)=f_!\ul{\Hom}_{\sO_{U,\bbox}}\Bigl(\Omega^m_{U/A}(\log|D|_U|),f^!A\otimes^{\L}_{\sO_{U,\bbox}}\sO_U(-D|_U+|D|_U|)\Bigr)\simeq g_!j_!j^*G.
\]
For the ambient morphism $j:(U,D|_U)\to(Y,D)$, the morphism $j^*\Omega^m_{Y/A}(\log|D|)(D-|D|)\to\Omega^m_{U/A}(\log|D|_U|)(D|_U-|D|_U|)$ used in Section~3 is the identity, and \eqref{solidmmorphism1} is the isomorphism of Lemma~\ref{compactprojection}. Hence, under the above identification, $j^m_+:\Solidm(U,D|_U)\to\Solidm(Y,D)$ is $g_!$ applied to the counit $j_!j^*G\to G$. Since $G\simeq F_m[d]$, Proposition~\ref{solidboundary} gives
\[
\cofib\Bigl(j^m_+:\Solidm(U,D|_U)\to\Solidm(Y,D)\Bigr)\simeq\Bigl(\bigl(g_{0*}\widehat F_{m,0}\bigr)\otimes^{\L}_{A_{0,\bbox}}A_{\bbox}\Bigr)[d].
\]
Finally, the functor $\ulMDMgm(A,\Q)\to D(A_{\bbox})$ in Theorem~\ref{mSolidthm} (7) is exact, and it sends the morphism $\ulM(U,D|_U)_{\Q}\to\ulM(Y,D)_{\Q}$ induced by $j$ to $j^m_+$. Therefore, it sends the cone to the above cofiber. This proves the claim.
\end{proof}

\begin{remark}
The isomorphism of cofibers in the proof of Theorem~\ref{solidboundarycone} holds over any field. The assumption that $A$ has characteristic $0$ is used only for the functor in Theorem~\ref{mSolidthm} (7).
\end{remark}

\bibliography{bib}

@misc{prosolid,
  author = {Ren, Fei},
  title = {Notes on Coherent Six-Functor Formalisms: Pro vs Solid},
  year = {2026},
  note = {arXiv:2506.21082v3},
  howpublished = {\url{https://arxiv.org/abs/2506.21082v3}}
}

@misc{cond,
  author = {Scholze, Peter},
  title = {Lectures on Condensed Mathematics},
  year = {2019},
  note = {All results joint with Dustin Clausen},
  howpublished = {\url{https://www.math.uni-bonn.de/people/scholze/Condensed.pdf}}
}

@misc{analyticstack,
  author = {Clausen, Dustin and Scholze, Peter},
  title = {Analytic Stacks},
  year = {2023--2024},
  note = {Course lectures},
  howpublished = {\url{https://people.mpim-bonn.mpg.de/scholze/AnalyticStacks.html}}
}

@article{miyazaki2020nisnevich,
  title={Nisnevich topology with modulus},
  author={Miyazaki, Hiroyasu},
  journal={Annals of K-theory},
  volume={5},
  number={3},
  pages={581--604},
  year={2020},
  publisher={Mathematical Sciences Publishers}
}

@misc{kelly2023hodgecohomologyramificationfiltration,
      title={Hodge cohomology with a ramification filtration, {II}}, 
      author={Shane Kelly and Hiroyasu Miyazaki},
      year={2023},
      eprint={2306.06864},
      archivePrefix={arXiv},
      primaryClass={math.AG},
      url={https://arxiv.org/abs/2306.06864},
      howpublished={\url{https://arxiv.org/abs/2306.06864}},
}

@article{kelly2023hodge,
  title={Hodge cohomology with a ramification filtration, {I}},
  author={Kelly, Shane and Miyazaki, Hiroyasu},
  journal={Mathematische Zeitschrift},
  volume={305},
  number={4},
  pages={70},
  year={2023},
  publisher={Springer}
}

@misc{ScholzeSix,
  author = {Scholze, Peter},
  title = {Six-Functor Formalisms},
  year = {2026},
  note = {arXiv:2510.26269v2},
  howpublished = {\url{https://arxiv.org/abs/2510.26269v2}}
}

@misc{Stacks,
  author = {{The Stacks Project Authors}},
  title = {The {Stacks Project}},
  year = {2026},
  howpublished = {\url{https://stacks.math.columbia.edu}}
}

@misc{CamargoSolid,
  author = {Rodr{\'i}guez Camargo, Juan Esteban},
  title = {Notes on Solid Geometry},
  year = {2026},
  note = {arXiv:2603.03012v1},
  howpublished = {\url{https://arxiv.org/abs/2603.03012v1}}
}
\bibliographystyle{unsrt}

%

%

\end{document}